# Decision Making: system lexicon.

Z. Alimbarashvili, V. Zhukovin, N. Chkhikvadze

**Abstract**: There is given the program realization system "lexicon" by the using lexicographical procedure.

## 1. Introduction

The lexicographical procedure widely spread in the decision making problems. It is used in the various sphere of the human action. for example:
1. The words are lexicographically regulated in the dictionary.
2. The base-ten system of the number recording is defined by lexicographical principle.
3. The militaries use the lexicographical principle in the problem of the objective allocation.
4. The American scientists experimentally proved, that the humans use the lexicographical principle for decision making.

## 2. Multicriteria lexicographical procedure.

Multicriteria decision making problem has been mathematical formalized with L. Zade in 1960. After that lexicographical procedure of choice has mathematical base for its description, development and generation new procedures.

Multicriteria problem of decision making is presented in vector interpretation in this way:

$$D = \langle X, K \rangle, \qquad (1)$$

where, $X$ - is finite set of competitive alternatives $x_i \subset X$, $i = 1 \div n$. Let introduce the set of ordered pairs $E$ with the elements $(x_i, x_k) \in E$, $x_i, x_k \in X$.

$$K = \{K_1(x_i), K_2(x_i), ..., K_j(x_i), ..., K_m(x_i)\} \qquad (2)$$

This is vector criterion of efficiency; where $K_j(x_i)$, $j = 1 \div m$, are scalar functions defined on the set $X$. We can suppose, that all $K_j(x_i)$ are type of "win". Besides they are defined with the content and although each of them has its scale of measuring the type of scale is identity.

Since we always can trade places the components at vector criterion $K$ of importance, we could say that $K_j(x_i)$ criterions are ordered by the relevance. The order concurs with the number's order, in other words most important criterion is $j = 1$, next $j = 2$ and so distances $j = m$. It is implied that this order is lexicographical. We explain this affirmation below. Now we formulate lexicographical structure on the set $E$. For it is necessary to define binary lexicographical preference relations $R_{lex}$. It defines in the following way. Let consider the decisions pairs $(x_i, x_l) \in E$.

$R_{lex}$: It is sad, that the decisions $x_i \in X$, are lexicographical preference than

$x_l \in X$, if one of following conditions fulfils:
1) $K_1(x_i) > K_1(x_l)$, or
2) $K_1(x_i) = K_1(x_l)$ and $K_2(x_i) > K_2(x_l)$, or  (3)

- - - - - - - - - - - - - - - - - -

j) $K_1(x_i) = K_1(x_l)$ and ...and $K_{j-1}(x_i) = K_{j-1}(x_l)$ and $K_j(x_i) > K_j(x_l)$ or

- - - - - - - - - - - - - - - - - - -

n) $K_1(x_i) = K_1(x_l)$ and ... and $K_{n-1}(x_i) = K_{n-1}(x_l)$ and $K_n(x_i) > K_n(x_l)$

The stop happens on the row where the inequality is fulfilled. If the inequality is not executed on $n$ row too, i.e. $K_n(x_i) = K_n(x_l)$, then this pair is lexicographically equivalent decision. Briefly this procedure is designating in that way:

$$x_i \overset{lex}{\underset{\sim}{\succ}} x_l \qquad (4)$$

Let explore characteristics of $R_{lex}$ In the first place it is **linked,** i.e. all decision pairs from $E$ are congruous with this relation. It is also **asymmetrical** and **transitive**.

Superiority relation $R_{lex}$ is equivalence too. Such superiority relation is linear order i.e. $R_{lex}$ is linear order. Now we shall apply all arsenal of result, theorems, which are available at present in the theories of decision making established us and other scientists, to research of the presented lexicographical procedure (3). And here the first. It is proved that the linear order has nonempty set Pareto (its kernel) and this set contains one decision, if it some they are equivalent. So procedure (3) always give out one decision

### 3. Lexicographical coefficients of importance.

Criterion $K_j(x_i)$ of efficiency are ordered on importance according to numbers: the less number the more important criteria for a choice. This fact can be fixed, having attributed criteria coefficients of importance: $\lambda_j \to K_j(x_i)$, $j = 1 \div m$. In this case in the theory of decision making linear convolution of vector criteria of efficiency $K$ is considered:

$$L(x_i) = \sum_{j=1}^{m} \lambda_j * K_j(x_i) \qquad (5)$$

This convolution possesses many very good properties which we shall consider later, in the same section. The unique requirement to importance coefficients are $\lambda_j \geq 0$ for all $j = 1 \div m$. Sometimes use also a condition $\sum_{j=1}^{m} \lambda_j = 1$, but it not definitely. With this the condition is executed:

if $L(x_i) \geq L(x_\ell)$ then $x_i \underset{\sim}{\succ} x_\ell$, and  (6)

on the contrary if $x_i \underset{\sim}{\succ} x_\ell$ then $L(x_i) \geq L(x_\ell)$.

Clear that if $\lambda_j$ are numbers and $K_j(x_i)$ scalar function, then $L(x_i)$ scalar function definite on $X$. The freedom of a choice $\lambda_j, j = 1 \div n$ in this case us is not arranged.

Self-wiled collection of important coefficients $\Lambda = \{\lambda_j\}_1^m$ disturbs lexicographical conditions representing by the formula (3). Whether there is a question is there such collection of importance coefficients which would not break, and kept a condition lexicographic (3)? Yes, it is and not only one but whole class. Therefore we designate lexicographical coefficients of importance with $\lambda_j(\ell ex)$, $j = 1 \div m$ and all collection of coefficients with

$$\{\lambda_j(\ell ex)\}_1^m = \Lambda(\ell ex) . \tag{7}$$

When linear convolution writes so for lexicographical choice (basic formula):

$$L_{\ell ex}(x_i) = \sum_{j=1}^m \lambda_j(\ell ex) * K_j(x_i) . \tag{8}$$

There are worked special procedures for the accounting of importance coefficients $\lambda_j(\ell ex)$, $j = 1 \div m$, by several scientists and by us too. These procedures do not coincide and use different scales of measuring. Now we describe basic property of linear convolution:

a) $L_{\ell ex}$ is a function of usefulness definite on the $X$. It assigns linear order according to rules:

$$L_{\ell ex}(x_i) \geq L_{\ell ex}(x_\ell) \leftrightarrow x_i \underset{\sim}{\succ}^{lex} x_\ell \tag{9}$$

b) This order is adjusted with the order presented by formula (3), this is:
   if $x_i \underset{\sim}{\succ} x_\ell$ correspondence with formula (3), that

$$L(x_i) \geq L(x_l) \tag{10}$$

and on the contrary. The formulas (8) and (9) are identical.

c) Any new presentation must be adjusted with (3) and it means that the presentation will be adjusted with (7) too.

d) Pareto set $X_\Pi(\ell ex)$ contains one decision (kernel of linear orders). This best decision is single and does not depend for the view of the presentation.

e) This best element will be find, if we decide simple optimization problem:

$$x^* = \max_{x_i \in X} L_{\ell ex}(x_i) \tag{11}$$

where $x^* \in X$ and it is best decision in lexicographical procedure of choice (3).

f) Now about scale of measure. All known us collection of lexicographic coefficients are presented by ordering scale:

$$Ш = \langle K_j, \Lambda_{\ell ex}, \Phi(K_j) \rangle , \tag{12}$$

with this $\Phi$ is the class of permitting transformation. It is monotone function (increasing or decreasing). Usually the decision $x^*$ is instable in the ordering scale i.e. the expression (10) does not invariant concerning monotone transformation of the coefficients. But this common rule breaks in lexicographical procedure The decision $x^*$ is stability concerning the transformation. The expression (10) is invariant in ordering scale. It is related to this max is searching with one criterion and which is defined by compare pair of decision.

g) Cleverly, lexicographical choice link to superiority degree: $Z(x_i, x_\ell) = j_0$
( it is number of line where is executed the choice).

Let analyze detail the criterions $K_j$. It is scalar function definite on $X$. Let consider these coefficients. Every criterion has own name (energy, length, expenditure), which represents any property of the estimating object $x_i \in X$. The name enables us the group of criterion $K_j, j = 1 \div m$, will be ordered with importance.

Every number criterion has the scale of estimating. There is example pic.1.

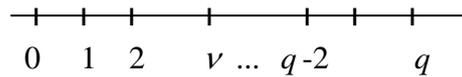

0  1  2   $v$ ...  q-2    q

pic.1.

Usually the scale is beginning from the zero, but it is possible beginning from either of meaning (for example $p$). Let consider the criterion of "win". It means that $v > v - 1$ where $v$ is indicator of scale rank. This dependence is transitive and linked. Now define scale diapason of criterion $K_j$.

$$d_j = \max v(j) - \min v(j), \qquad (13)$$

if the scale is beginning from zero $d_j = \max v(j) = q$, where $q$ is quantity of scale rank.

Let define the conditions of lexicographical group criterions $K_j, j = 1 \div m$:

**Affirmation 1.**

Assume the importance of criterion is increasing from the right to the left and its number of position from the left to the right (formulae 3), then if the condition $\min v(j) \succ \max v(j+1)$ is fulfilled for all $j = 1 \div m - 1$, then criterions group $K_j, j = 1 \div m$ is lexicographic ordered.

Notice:
1. Zero takes not in part of the minimum defining (only in the affirmation)
2. The mark $\succ$ signifies the "importance"
3. This dependence is linked and transitivity.

Now we can assign the formula for the accounting lexicographical coefficients of importance:

$$\lambda_j(\ell ex) = d_j^{j-1}, \quad j = 1 \div m \qquad (14)$$

May be the scales of the criterion are diverse (for example: the speed m/h, energy k/h). In this case the formulae (3) works, but $L(\ell ex)$ loses the sense. One of the method the forming uniformity system of the criterion is the rationing. Instead of $K_j$, include following criterions:

$$\gamma_j(x_i) = \frac{K_j(x_i)}{\max K_j} * a \quad \text{where } a > 0 \qquad (15)$$

This gives the freedom for the variation. ( for example $a = 10$ scale of Saaty).

**Uniformity of criteria.** Let consider two example of lexicographical structure for the presenting uniformity criterions.

a) Decimal system of number presenting:

$15530 = 1*10^4 + 5*10^3 + 3*10^2 + 0*10^0$

The position is corresponding with the criterion. The numbering of the position is from the right to the left. The importance grows up same direction too. There is same the scale in each position. This is the uniformity of the criterion (position)

b) Lexicon.

The direction of importance increasing is from the right to the left and position numbering on the contrary from the left to the right. The number of the position is unlimited. It is bad for the dictionary. For this we assign maximal possible length of the word ( 20 letter). One letter contains in each position ("gap" is placed in the beginning of the alphabet or in end). The letters are ordered by the preference (from A to Z). There is uniformity structure too.

**Notice:** Each word is ascribed the number then we order these numbers. The accounting lexicographical coefficients of importance are for Georgian alphabet:

$$\lambda_j(\ell ex) = 33^{j-1} \text{ for all } j.$$

### 4. System lexicon.

Now because our article is about the lexicon (Georgian-English, English-Georgian, Russian-Georgian, Georgian- Russian i.e.). we consider one of the example for the computer.

There are some restriction limitation which we permitted during the programming.
1. The number of positions is 20 .
2. The base of preference we take 6 (instead of 33 in Georgian and 26 in English).
3. First position of preference is $\lambda = 6^{30}$ (top limit of the memory ).
4. The preference is increasing from the left to the right.
5. For the unique of each position $\lambda^{j+1} - \lambda^j = 6^2$ .
6. The diapason of preference difference is from $\lambda = 6^{30}$ - to $\lambda = 6^{-8}$, because number of position is 20.

The program part of the system "lexicon" is formed in system **Delphi .** system "lexicon" composes two programs and they work dialogic regime.

First of the programs "main" search the font of the words, the context and the corresponding data base. After all it pasts the ruling second program.

Second program is "lexicon". There are the commentaries which realize recess operation. There are considered:
1. Get the word and context.
2. **The word and context are effaced if they are wrong.**
3. **Special code weight and places corresponding to its place in data base.**
4. **The lexicon may be carried to EXCEL and it gives possibility the publication lexicon for wide customer.**